\documentclass[12pt]{article}
\usepackage{amsmath,amssymb,amsthm}
\def\abs#1{\left|#1\right|}
\def\Cnt#1{{\cal C}^{#1}}
\def\comp{\circ}
\def\C{\mathbb C}
\def\csub{\subset\subset}
\def\eps{\varepsilon}
\def\inv#1{{#1}^{-1}}
\def\meas{\mu}
\def\Lor#1#2{SO^+(1,#1,#2)}
\def\N{\mathbb N}
\def\nsR{{\widetilde{\mathbb R}}_c}
\def\Ord#1{O#1}
\def\Q{\mathbb Q}
\def\R{\mathbb R}
\def\Rot#1#2{SO(#1, #2)}
\def\rot{\mathcal R}
\def\sign{\sigma}
\def\Trans#1{T_{#1}}
\def\transp#1{{#1}^t}
\def\Z{\mathbb Z}

\def\Gen{\mathcal G}
\def\Mod{{\mathcal E}_M}
\def\Null{\mathcal N}

\def\ns{\mathrm{ns}}
\def\ster#1{{{}^* \mskip-1mu #1}}

\newtheorem*{cor}{Corollary}
\newtheorem{df}{Definition}
\newtheorem{thm}{Theorem}
\newtheorem{lemma}[thm]{Lemma}
\newtheorem*{rem}{Remark}

\begin{document}
\title{Group invariant Colombeau generalized functions}

\author{Hans Vernaeve\footnote{Supported by research grant M949 of the Austrian Science Foundation (FWF)}\\
Institut f\"ur Grundlagen der Bauingenieurwissenschaften\\
Technikerstra\ss e 13\\
A-6020 Innsbruck}
%E-mail: {\tt Hans.Vernaeve@uibk.ac.at}}

\date{}
\maketitle

\begin{abstract}
Colombeau generalized functions invariant under smooth (additive) one-parameter groups are characterized.
This characterization is applied to generalized functions invariant under orthogonal groups of arbitrary signature, such as groups of rotations or the Lorentz group. Further, a one-dimensional Colombeau generalized function with two (real) periods is shown to be a generalized constant, when the ratio of the periods is an algebraic nonrational number. Finally, a nonstandard Colombeau generalized function invariant under standard translations is shown to be constant.
\end{abstract}



\section{Introduction}
This paper is related to a series of papers on group invariant generalized functions that appeared during the last years (\cite{Kunzinger06}, \cite{Kunzinger06b}, \cite{Kunzinger97}, \cite{Kunzinger00}, \cite{Ober04}). In particular, this paper focuses on a type of questions that remained an open problem for several years: if a generalized function is invariant under all non-generalized transformations of a generalized transformation group, is it then invariant under the whole (generalized) group? Only recently, the key case of a translation group was solved in the affirmative in~\cite{PSV}. In \cite{Kunzinger06}, it is already shown that this result can be applied to solve the question for the group of rotations. In this paper, we show that it can be applied to more general groups as well. We prove a general result on invariance under smooth one-parameter groups (section~\ref{section oneparam}) and indicate how it can be applied to invariance under various matrix groups (sections~\ref{section rotation}, \ref{section Lorentz}). In the case of rotations, the same characterization as in \cite{Kunzinger06} is obtained. This development mirrors results on group invariant Schwartz distributions, following the work by Schwartz~\cite{Schwartz54} (see~\cite{Berest}, \cite{Ibragimov}, \cite{Methee54}, \cite{Szmydt}, \cite{Tengstrand60} and references therein).\\
In section~\ref{section Lobster}, we revisit the case of the translation groups. We give two more proofs of this theorem. Doing so, the following new results are obtained: a one-dimensional Colombeau generalized function with two (real) periods is a generalized constant, when the ratio of the periods is an algebraic nonrational number; a nonstandard Colombeau generalized function invariant under standard translations is constant.

\section{Preliminaries}
We work in the (so called {\it special}) Colombeau algebra $\Gen(\R^d)$ of generalized
functions on $\R^d$ ($d\in\N$), defined as follows \cite{Col, GKOS}.\\
Denote by ${\mathcal E}(\R^d)$ the algebra of all nets $(u_\eps)_{\eps\in(0,1)}$
of $\Cnt{\infty}$-functions $\R^d\to\C$.\\
Then $\Gen(\R^d)=\Mod(\R^d)/\Null(\R^d)$, where
\begin{align*}
\Mod(\R^d) =\,&\big\{(u_\eps)_\eps \in{\mathcal E}(\R^d): (\forall
K\csub\R^d) (\forall\alpha\in\N^d) (\exists b\in\R)\\
&\big(
\sup_{x\in K}\abs{\partial^\alpha u_\eps(x)}=\Ord(\eps^b),\text{ as }\eps\to 0
\big)\big\}\\
\Null(\R^d) =\,&\big\{(u_\eps)_\eps \in{\mathcal E}(\R^d): (\forall
K\csub\R^d) (\forall\alpha\in\N^d) (\forall b\in\R)\\
&\big(
\sup_{x\in K}\abs{\partial^\alpha u_\eps(x)}=\Ord(\eps^b),\text{ as }\eps\to 0
\big)\big\}.
\end{align*}
We recall a result about the composition of generalized functions.
\begin{df}
An element $(f_\eps)_\eps = (f_{1,\eps},\dots, f_{d,\eps})$ $\in\Mod(\R^d)^d$ is
called c-bounded if
\[(\forall K\csub\R^d)(\exists K'\csub\R^d)(\exists\eps_0>0)
(\forall\eps<\eps_0)(f_\eps(K)\subseteq K').
\]
An element of $\Gen(\R^d)^d$ is called c-bounded if it possesses a c-bounded
representative.
\end{df}

\begin{lemma}
1. Let $f\in \Gen(\R^d)^d$ be c-bounded and $g\in\Gen(\R^d)$. Then the composition
$g\comp f$ defined on representatives by means of
\[(g\comp f)_\eps = g_\eps\comp f_\eps\]
is a well-defined generalized function in $\Gen(\R^d)$.\\
2. Let $f$, $g$ $\in\Gen(\R^d)^d$ be c-bounded. Then $g\comp f$ (similarly
defined on representatives) is a well-defined c-bounded generalized function in
$\Gen(\R^d)^d$.
\end{lemma}
\begin{proof}
1. See~\cite[section 1.2.1]{GKOS}.\\
2. The well-definedness follows directly from the first part; the c-boundedness
of $g\comp f$ follows directly from $(g\comp f)_\eps(K)\subseteq g_\eps(f_\eps(K))$ and
the c-boundedness of $f$ and $g$.
\end{proof}

In the last section, we will also work in the algebra of nonstandard Colombeau generalized functions ${}^\rho{\mathcal E}(\R^d)$, defined as follows \cite{Ober88, Todorov88}.\\
Let $\rho\in\ster\R$ be a fixed positive infinitesimal.\\
For $x,y\in\ster\C^d$, we write $x\approxeq y$ iff $\abs{x-y}\le\rho^n$,
$\forall n\in\N$ ($x-y$ is then called a iota, or negligible).\\
For $x\in\ster\C^d$, we write $x\in\ster\C_M$ iff $\abs{x}\le 1/\rho^n$ for
some $n\in\N$ ($x$ is then called moderate).\\
We denote by $\ns(\ster\R^d)$ the set of near-standard (=finite) elements of $\ster\R^d$.\\
%The field $^\rho\C:=\ster\C/\approxeq$ is called the field of (complex) asymptotic numbers.\\
Then ${}^\rho{\mathcal E}(\R^d)=\Mod(\R^d)/\Null(\R^d)$, where
\begin{align*}
\Mod(\R^d)&=\{u \in\ster{\Cnt\infty}(\R^d): (\forall x\in\ns(\ster\R^d)) (\forall \alpha\in\N^n)
(\partial^\alpha u(x)\in\ster\C_M)\}\\
\Null(\R^d)&=\{u \in\ster{\Cnt\infty}(\R^d): (\forall x\in\ns(\ster\R^d))(\forall \alpha\in\N^n)
(\partial^\alpha u(x)\approxeq 0)\}.
\end{align*}

\section{Invariance under one-parameter groups}\label{section oneparam}
\begin{thm}\label{one-param}
Let $g_\theta$ be a one-parameter (additive) group action on $\R^d$, i.e.
\begin{enumerate}
\item $g_\theta$ is a bijection, $\forall \theta\in\R$
\item $g_{\theta_1+\theta_2} =  g_{\theta_1}\comp g_{\theta_2}$, $\forall \theta_1, \theta_2\in\R$
\item $g_0$ is the identity mapping on $\R^d$
\item $g_{-\theta}=\inv{(g_\theta)}$, $\forall \theta\in\R$.
\end{enumerate}
Suppose further that the map $(\theta, x)\mapsto g_\theta(x)$ is
$\Cnt{\infty}$.\\
Let $f\in \Gen(\R^d)$. Then for each $\theta \in\nsR$, $f\circ g_\theta$ is well-defined.\\
If $f\circ g_\theta = f$, $\forall \theta \in\R$, then $f\circ g_\theta = f$, $\forall \theta \in\nsR$, where
$\nsR$ is the ring of Colombeau generalized numbers with bounded representative.
\end{thm}
\begin{proof}
First, notice that the $\Cnt\infty$-character of $(\theta,x)\mapsto g_\theta$
implies that
\begin{equation}
\label{c-bounded}
(\forall k\in\N)(\forall\alpha\in\N^d)(\forall R\in\R)
\Big(\sup_{\abs{\theta}\le R}\sup_{\abs{x}\le R}
\abs{\frac{\partial^k}{\partial \theta^k}\partial^\alpha g_\theta(x)}<+\infty
\Big)
\end{equation}
so in particular, $g_\theta$ is a c-bounded generalized function, for each $\theta\in\nsR$ and the
composition is well-defined.\\
Fix $a = (a_\eps)_\eps\in\nsR^d$.
Define $h\in\Gen(\R)$ by means of
$h(\theta) = f(g_\theta(a))$, i.e., by a definition on representatives,
\[h_\eps(\theta) = f_\eps(g_\theta(a_\eps)).\]
By the c-boundedness of $g_\theta$, one sees that this
definition is independent of the representative of $f$. Further, by
equation~(\ref{c-bounded}), one also sees that $(h_\eps)_\eps\in\Mod(\R)$.\\
Now let $\Trans{t}$ be the translation-operator $x\mapsto x+t$ on $\R$. Then
\[
(h\comp\Trans{t})(\theta) = f(g_{\theta+t}(a)) = (f\comp g_t) (g_\theta(a)),
\]
so by the hypothesis on $f$
\[
(h\comp\Trans{t})(\theta) = f(g_\theta(a)) = h(\theta).
\]
%[To be explicit, one should write the previous equality on representatives to see that
%the equality $f\comp g_t = f$ on $\Gen(\R^d)$ yields the equality
%$h\comp\Trans{t} = h$ in $\Gen(\R)$.]\\
So $h$ is a generalized constant (see section~\ref{section Lobster}). I.e., for each $C\in\R^+$,
\begin{equation}
(\forall p \in\N)(\exists \eps_0\in\R^+)(\forall \eps<\eps_0)
\Big(
\sup_{\substack{\theta\in\R\\\abs{\theta}<C}}
\abs{h_\eps(\theta)-h_\eps(0)}\le \eps^p
\Big).
\label{to-contradict}
\end{equation}
We conclude that, for each $R\in\R^+$,
\[
(\forall p \in\N)(\exists \eps_0\in\R^+)(\forall \eps<\eps_0)
\Big(
\sup_{\substack{x\in\R^d\\\abs{x}<R}}
\underbrace{
\sup_{\substack{\theta\in\R\\\abs{\theta}<C}}
\abs{f_\eps(g_{\theta}(x))-f_\eps(x)}}
_{= F_\eps(x)} \le \eps^p
\Big),
\]
since supposing the contrary, one can construct, for some $p\in\N$, a decreasing
sequence $(\eps_n)$ with $\eps_n < 1/n$ such that there exists
$a_{\eps_n}\in\R^d$ with
$\abs{a_{\eps_n}}<R$ and $F_{\eps_n}(a_{\eps_n}) > \eps_n^p$. Extending $(a_{\eps_n})_n$
arbitrarily into $(a_\eps)_\eps$ (with $\abs{a_\eps}<R$, $\forall\eps$) yields
then $a\in\nsR^d$ for which equation~(\ref{to-contradict}) is false.\\
Now let $\bar\theta\in\nsR$. Then $(f\comp g_{\bar\theta})_\eps (x) =
f_\eps(g_{\theta_\eps}(x))$ with, for some $C\in\R^+$,
$\abs{\theta_\eps}\le C$, $\forall\eps$. We have in particular that
\[
(\forall p \in\N)(\exists \eps_0\in\R^+)(\forall \eps<\eps_0)
\Big(
\sup_{\substack{x\in\R^d\\\abs{x}<R}}
\abs{f_\eps(g_{\theta_\eps}(x))-f_\eps(x)} \le \eps^p
\Big),
\]
i.e.\ (as $R$ can be chosen arbitrarily), $f\comp g_{\bar\theta} = f$.
\end{proof}
\begin{cor}
%\begin{enumerate}
1. Let $\alpha$ be a fixed plane through the origin in $\R^d$ (i.e., a
2-dimensional subspace of $\R^d$).
Let $g_\theta$ be the rotation in $\R^d$ over the angle $\theta$ in $\alpha$
(with a chosen orientation). Then $g_\theta$ is a one-parameter group action.
Explicitly, after a change of basis, $g_\theta$ is the linear transformation
with matrix
\[
\begin{pmatrix}
\cos\theta &-\sin\theta & 0 &\cdots & 0\\
\sin\theta &\cos\theta & 0 & \cdots & 0\\
0& 0\\
\vdots&\vdots&&I\\
0&0
\end{pmatrix}\,.
\]
It follows that $(\theta,x)\mapsto g_\theta(x)$ is a $\Cnt{\infty}$-mapping.
So, by the previous theorem, if $f\in\Gen(\R^d)$ is invariant under all rotations
$\in \Rot{d}{\R}$ in the plane $\alpha$, then $f$ is also invariant under all rotations
$\in \Rot{d}{\nsR}$ in the plane $\alpha$.\\
2. Let $\sigma_{i,j,\theta}$ be the linear transformation on $\R^d$ defined on coordinates by
\begin{multline*}
\sigma_{i,j,\theta}(x_1,\dots, x_d)= (x_1, \dots, x_{i-1}, x_i \cosh\theta + x_j \sinh\theta,\\
x_{i+1},\dots, x_{j-1}, x_i \sinh\theta + x_j \cosh\theta, x_{j+1}, \dots, x_d).
\end{multline*}
%\[
%\begin{pmatrix}
%\cosh\theta &\sinh\theta & 0 &\cdots & 0\\
%\sinh\theta &\cosh\theta & 0 & \cdots & 0\\
%0& 0\\
%\vdots&\vdots&&I\\
%0&0
%\end{pmatrix}\,.
%\]
Then also $\sigma_{i,j,\theta}$ is a one-parameter group action and $(\theta,x)\mapsto
\sigma_{i,j,\theta}(x)$ is a $\Cnt{\infty}$-mapping; so the previous theorem applies.
%\end{enumerate}
\end{cor}

\section{Invariance under groups of rotations}\label{section rotation}
In order to lift the theorem about one-parameter group actions to certain higher
dimensional group actions, we use the following elementary lemma.
\begin{lemma}
Let $(e_1$, \dots, $e_d)$ be the standard basis in $\R^d$.
Denote by $\rot_{i,j,\theta}$ (with $i < j$) the rotation over the angle $\theta$ in the
vectorplane spanned by the vectors $e_i$, $e_j$ (with a chosen orientation).\\
There exist $i_1$, \dots, $i_{\binom{d}{2}}$ and $j_1$, \dots, $j_{\binom{d}{2}}$
such that every $\rot\in \Rot{d}{\R}$ equals
\[\rot_{i_1,j_1,\theta_1}\comp \dots\comp
\rot_{i_{\binom{d}{2}},j_{\binom{d}{2}},\theta_{\binom{d}{2}}},\]
for some $\theta_1$, \dots, $\theta_{\binom{d}{2}}$ $\in\R$.\\
Similarly, every $\rot\in\Rot{d}{\nsR}$ equals
\[\rot_{i_1,j_1,\theta_1}\comp \dots\comp
\rot_{i_{\binom{d}{2}},j_{\binom{d}{2}},\theta_{\binom{d}{2}}},\]
for some $\theta_1$, \dots, $\theta_{\binom{d}{2}}$ $\in\nsR$.
\end{lemma}
\begin{proof}
Let $\rot\in\Rot{d}{\R}$ and call $x=\rot(e_d)$. Apply $\rot_{1,d,\theta_1}$ with $\theta_1$ such that, after
rotating, $x_1=0$; apply $\rot_{2,d,\theta_2}$ with $\theta_2$ such that, after
rotating, $x_1=x_2=0$; \dots
So we can find $\theta_1$, \dots, $\theta_{d-1}$ such that
$\rot_{d-1,d,\theta_{d-1}}\comp\cdots\comp\rot_{2,d,\theta_2}\comp\rot_{1,d,\theta_1}\comp\rot$
is a rotation that fixes $e_d$, i.e., a rotation in $\{x\in\R^d: x_d = 0\}\cong\R^{d-1}$.
By induction,
\[
\rot_{1,2,\theta_{\binom{d}{2}}}\comp(\rot_{2,3,\theta_{\binom{d}{2}
-1}}\comp\rot_{1,3,\theta_{{\binom{d}{2}}-2}})\comp\cdots\comp
(\rot_{d-1,d,\theta_{d-1}}\comp\cdots\comp\rot_{2,d,\theta_2}\comp\rot_{1,d,\theta_1})\comp\rot
\]
is the identity (for suitable $\theta_1$, \dots, $\theta_{\binom{d}{2}}$). I.e., $\rot$ equals
\[
\rot_{1,d,-\theta_1}\comp\rot_{2,d,-\theta_2}\comp\cdots\comp\rot_{d-1,d,-\theta_{d-1}}
\comp\cdots\comp\rot_{1,3,-\theta_{{\binom{d}{2}}-2}}
\comp\rot_{2,3,-\theta_{{\binom{d}{2}}-1}}
\comp\rot_{1,2,-\theta_{\binom{d}{2}}},
\]
so $\rot$ has the required form.\\
Now let $\rot\in\Rot{d}{\nsR}$. It has a representative $(\rot_\eps)_\eps$, with
$\rot_\eps\in\Rot{d}{\R}$, $\forall\eps$. So $\rot_\eps$ equals, for suitable
$\theta_{i,\eps}\in[0,2\pi]$,
\[
\rot_{1,d,-\theta_{1,\eps}}\comp\rot_{2,d,-\theta_{2,\eps}}\comp\cdots
\comp\rot_{1,3,-\theta_{{\binom{d}{2}}-2,\eps}}
\comp\rot_{2,3,-\theta_{{\binom{d}{2}}-1,\eps}}
\comp\rot_{1,2,-\theta_{{\binom{d}{2}},\eps}}.
\]
Since each $\rot_{i,j,\bar\theta}$ (with $\bar\theta\in\nsR$) is c-bounded, all
compositions of generalized functions are well-defined, so we conclude that
$\rot$ has the required form.
\end{proof}

We obtain the same answer as in~\cite{Kunzinger06} to the open question posed in~\cite{Ober04}.
\begin{thm}\label{full-rot}
Let $f\in \Gen(\R^d)$, $f\comp \rot = f$, $\forall \rot \in\Rot{d}{\R}$.
Then $f\comp \rot = f$, $\forall \rot \in\Rot{d}{\nsR}$.
\end{thm}
\begin{proof}
By corollary~1 to thm.~\ref{one-param}, $f$ is invariant under all rotations
$\in \Rot{d}{\nsR}$ in a fixed plane. In particular, $f=f\comp \rot_{i,j,\theta}$,
for all rotations $\rot_{i,j,\theta}$ as defined in the previous lemma
($\theta\in\nsR$). Now let $\rot\in\Rot{d}{\nsR}$ arbitrarily. Since $\rot$
equals
$\rot_{i_1,j_1,\theta_1}$ $\comp \cdots\comp$
$\rot_{i_{\binom{d}{2}},j_{\binom{d}{2}},\theta_{\binom{d}{2}}}$,
for some $\theta_1$, \dots, $\theta_{\binom{d}{2}}$ $\in\nsR$, we conclude (since
all compositions are well-defined because of c-boundedness) that
\begin{align*}
f\comp \rot
&=(f\comp\rot_{i_1,j_1,\theta_1})\comp\rot_{i_2,j_2,\theta_2}\comp\cdots\comp
\rot_{i_{\binom{d}{2}},j_{\binom{d}{2}},\theta_{\binom{d}{2}}}\\
&= f\comp\rot_{i_2,j_2,\theta_2}\comp\cdots\comp\rot_{i_{\binom{d}{2}},j_{\binom{d}{2}},\theta_{\binom{d}{2}}}
= \cdots = f.
\end{align*}
\end{proof}

\begin{rem}
A similar result holds for the full orthogonal group $O(d,\R)$ (i.e., the group of all
linear transformations on $\R^d$ preserving the usual inner product), since
any element from this group can be represented as an element
of $\Rot{d}{\R}$, possibly composed with one fixed
orientation-inverting orthogonal transformation.
\end{rem}

\section{Invariance under Lorentz-transformations}\label{section Lorentz}
Let $\Lor{d}{\R}$ be the group of all proper, orthochronous Lorentz transformations
in $\R^{d+1}$, i.e., the group of all linear transformations of $\R^{d+1}$
preserving the form $t^2-\abs{x}^2$, $(t,x)\in\R^{d+1}$, as well as orientation,
and the direction of time (the last two assertions mean that for the matrix
$A=(a_{ij})$ of the linear transformation, $\det A= + 1$ and $a_{00} > 0$).\\
Then the following is true (see~\cite[Appendix]{Szmydt}).
\begin{lemma}
Every $g\in\Lor{d}{\R}$ equals $\rot_1 \comp \sigma_{1,2,\theta}\comp \rot_2$, where
$\rot_1$,$\rot_2$ are rotations in $\Rot{d}{\R}$ (which keep the time variable $t$
invariant), and $\sigma_{1,2,\theta}$ ($\theta\in\R$) is defined as in the corollary to
theorem~\ref{one-param} (here, the first variable means the time variable $t$).
\end{lemma}
Again, we can transfer this result to the analogous result about
$\Lor{d}{\nsR}$, the group of (generalized) Lorentz transformations with
coefficients in $\nsR$:\\
Every $g\in\Lor{d}{\nsR}$ equals $\rot_1 \comp \sigma_{1,2,\theta}
\comp \rot_2$, for some $\rot_1$, $\rot_2$ $\in\Rot{d}{\nsR}$ and
$\theta\in\nsR$.
\begin{thm}
Let $f\in\Gen(\R^{d+1})$, $f\comp g = f$, $\forall g\in\Lor{d}{\R}$.
Then $f\comp g = f$, $\forall g\in\Lor{d}{\nsR}$.
\end{thm}
\begin{proof}
Since $f$ is invariant under all rotations in $\Rot{d}{\R}$ and under all
$\sigma_{1,2,\theta}$ ($\theta\in\R$), the previous results show that $f$ is also
invariant under all rotations in $\Rot{d}{\nsR}$ and under all $\sigma_{1,2,\theta}$
($\theta\in\nsR$).
Then $f$ is also invariant under any $\rot_1\comp\sigma_{1,2,\theta}\comp\rot_2$
($\rot_1,\rot_2\in\Rot{d}{\nsR}$; $\theta\in\nsR$).
\end{proof}
\begin{rem}
1. A similar result holds for the full Lorentz group (i.e., the group of all
nonsingular linear transformations preserving the form $t^2-\abs{x}^2$), since
any element from this group can be represented as an element
of $\Lor{d}{\R}$, possibly composed with one fixed time-inverting
Lorentz-transformation and possibly composed with one fixed
orientation-inverting Lorentz-transformation.\\
2. More generally, let $B$ be a nondegenerate real, symmetric bilinear form on $\R^d$. Let $A$ be a linear transformation that leaves $B$ invariant, i.e., in matrix notation, $\transp{A}BA=B$. Let $f\in\Gen(\R^d)$ be invariant under $A$, i.e., $f\comp A = f$. Then, for any invertible linear transformation $C$, $g=f\comp C\in\Gen(\R^d)$ is invariant under $\inv{C} A C$. So, $f$ is invariant under any linear transformation that leaves $B$ invariant iff $g$ is invariant under any transformation that leaves $\transp{C} B C$ invariant. So we may suppose that $B$ is reduced into its standard form
\[B(x,y)=x_1 y_1 +\cdots + x_p y_p - (x_{p+1}y_{p+1} + \cdots + x_d y_d), \quad p\le d.\]
Let $O(p,q,\R)$ ($q=d-p$) be the group of all linear transformations on $\R^d$ that leave $B$ invariant. Since a linear transformation in $O(p,q,\R)$ can always be written as $A_1\circ A_2\circ \sigma_{i,j,\theta}$, where $A_1\in O(p,\R)$ leaves $x_{p+1}$, \dots, $x_d$ invariant, $A_2\in O(q,\R)$ leaves $x_1$, \dots, $x_p$ invariant, $i\le p$ and $j>p$ (see e.g.\ \cite{Tengstrand60}), a similar reduction can be applied and we obtain that $f\in\Gen(\R^d)$ is invariant under $O(p,q,\R)$ iff $f$ is invariant under $O(p,q,\nsR)$.
%explicitly, the reduction of rotations in the elementary lemma in the previous paragraph can be carried out in this case as well: let $A(e_d)=x$; then reduce to a vector $(x_1, \dots, x_p, 0, \dots, 0, x_d)$ by means of rotations; as the bilinear form $B$ is preserved, $x_1^2 + \cdots + x_p^2 + 1 = x_d^2$. In particular, $x_d^2 > x_1^2$, which can be reduced (looking at the hyperbolic orbits of $\sigma_{i,j,\theta}$) by applying some $\sigma_{1,d,\theta}$ to $(0, x_2, \dots, x_p, 0, \dots, 0, x_d)$. Similarly, we can reduce to $(0,\cdots, 0, x_d)$ with $x_d=\pm 1$. Using a reflection, we can reduce to $(0, \cdots, 0, 1)$. Again (as the matrix of $A$ leaves $B$ invariant), this means that we are essentially reduced to the $d-1$-dimensional case.
\end{rem}

\section{Translation invariance}\label{section Lobster}
We revisit the theorem on translation-invariant generalized functions that was
used in the proof of theorem~\ref{one-param}. It was first proved in~\cite{PSV}.
We give two original proofs.
\begin{thm}
Let $u\in\Gen(\R^d)$.
Suppose that $u$ is invariant under translations $x\mapsto x + h$, $\forall h\in\R^d$. I.e.,
$u(x+h)=u(x)$ holds as an equality in $\Gen(\R^d)$, $\forall h\in\R^d$. Then $u$ is constant in $\Gen(\R^d)$.
\end{thm}
\begin{proof}
We suppose that $u$ is not constant, i.e., for a representative $(u_\eps)_\eps$ of $u$, $(u_\eps(x) - u_\eps(0))_\eps\notin\Null(\R^d)$.
So there exist $K\csub\R^d$, $N\in\N$ and a decreasing sequence
$(\eps_n)_{n\in\N}$, with $\lim_{n\to\infty}\eps_n = 0$, on which
\[\sup_{x\in K} \abs{u_{\eps_n}(x)-u_{\eps_n}(0)} > \eps_n^N.\]
Now let $f_n:= u_{\eps_n}-u_{\eps_n}(0)$. So there exists a sequence $(a_n)_{n\in\N}$,
$a_n\in K$ such that
\[\abs{f_n(a_n)} > \eps_n^N,\]
for each $n$. We have $\|a_n\|<C$, for some $C\in\R^+$.\\
Now let
\[A_n:= \Big\{x\in\R^d: \Big(\abs{f_n(x)}< \frac{\eps_n^N}{3}\Big)\Big\},\quad
B_n := \bigcap_{m\ge n} A_m.\]
As $f_n(0)=0$, $\forall n$, $0\in B_0$. Let $x\in\R^d$. Because of the
translation-invariance, also $x\in B_n$, for some $n$.
Similarly, we define
\[C_n:= \Big\{x\in\R^d: \Big(\abs{f_n(x + a_n)-f_n(a_n)}<
\frac{\eps_n^N}{3}\Big)\Big\},\quad
D_n := \bigcap_{m\ge n} C_m.\]
Let $x\in\R^d$. Again, $x\in D_n$, for some $n$.\\
Both $(B_n)$, $(D_n)$ are increasing sequences of measurable subsets of $\R^d$ and
$\bigcup B_n=\bigcup D_n = \R^d$. We denote the Lebesgue measure by $\meas$
and the open ball with center $x\in\R^d$ and radius $r\in\R^+$ by $B(x,r)$.
Now \[\meas(B(0,C)\setminus B_n)\to 0, \, \meas(B(0, 2C)\setminus D_n)\to 0\] as $n\to\infty$.
Since $B_n \subseteq A_n$, $D_n\subseteq C_n$, also
\[\meas(B(0,C)\setminus A_n)\to 0,\, \meas(B(0, 2C)\setminus C_n)\to 0\] as $n\to\infty$.
Finally, let
\[E_n := \Big\{x\in\R^d: \Big(\abs{f_n(x)-f_n(a_n)}< \frac{\eps_n^N}{3}\Big)\Big\} = C_n + a_n.\]
Then
\[\meas (B(0,C)\setminus E_n) \le \meas (B(a_n, 2C)\setminus E_n)
= \meas (B(0, 2C)\setminus C_n) \to 0,\]
since $\|a_n\|< C$. So $\meas (B(0,C)\setminus (A_n\cap E_n))\to 0$.
In particular, $A_n$ and $E_n$ have a non-empty intersection as soon as
$n$ is large enough. Clearly, this is impossible.
\end{proof}

The second proof is for the one-dimensional case (the more-dimensional case can then be obtained e.g.\ in a similar way as we obtained theorem~\ref{full-rot} by applying theorem~\ref{one-param}, or as in the corollary to theorem~\ref{ns-transinvar}) and uses weaker hypotheses on the generalized function. Moreover, both the previous proof and the proof given in~\cite{PSV} cannot (at least not a priori) be generalized to the nonstandard case (the ultrafilter destroys the argument). Although the proof of the nonstandard theorem is a little more conceptual, the proof of the standard version doesn't make use of nonstandard analysis.\\
We recall two number theoretic theorems.

\begin{thm}[Dirichlet's approximation theorem]\label{Dirichlet}
Let $\alpha\in\R^+\setminus\Q$. Then
\[
(\forall N\in\N\setminus\{0\})(\exists k,l\in\N)
\left(0<l\le N \quad\& \quad\abs{k-l\alpha}\le\frac{1}{N}\right).
\]
\end{thm}
\begin{proof}
See~\cite{Apostol}.
\end{proof}

\begin{thm}[Liouville's approximation theorem]\label{Liouville}
Let $\alpha\in\R^+\setminus\Q$ be an algebraic number of degree $n$.
Then there exists $c\in\R^+$  such that for each $k$, $l$ $\in\N$ ($l\ne 0$),
$\abs{\alpha - \frac{k}{l}}\ge\frac{c}{l^n}$.
\end{thm}
\begin{proof}
See~\cite{Apostol}.
\end{proof}

For our application, we will use the following corollary:
\begin{cor}
Let $\alpha\in\R^+\setminus\Q$ be an algebraic number.
Then there exists $M\in\N$ such that
\[
(\forall R\in(2,+\infty))(\exists k,l\in\N)
\left(l\le R \quad\& \quad \frac{1}{R^M}\le\abs{k-l\alpha}\le\frac{2}{R}\right).
\]
\end{cor}
\begin{proof}
Let $R\in\R$, $R\ge 2$. Let $N\in\N$ such that $R-1\le N \le R$. Then by
Dirichlet's approximation theorem, there exist $k,l\in\N$ such that $0<l\le R$
and $\abs{k-l\alpha}\le\frac{1}{R-1}\le\frac{2}{R}$.
By Liouville's approximation theorem, there exist $c\in\R^+$ and $n\in\N$
($n\ge 2$, both
depending on $\alpha$ only) such that $\abs{k-l\alpha}\ge \frac{c}{l^{n-1}}\ge
\frac{c}{R^{n-1}}\ge \frac{1}{R^M}$, for a good choice of $M$ (depending on $c$ and $n$,
hence on $\alpha$ only).
\end{proof}

We still need the following elementary lemma.
\begin{lemma}\label{almost-periodic}
Let $f$: $\R\to\C$ be almost-periodic on an interval $[a,b]$
with periods $h_1$,$h_2$ $\in\R^+$ and tolerance $\eps\in\R^+$, i.e., for
$i=1,2$
\[
(\forall x\in[a,b])(x+h_i\in [a,b]\implies \abs{f(x+h_i)-f(x)}\le\eps).
\]
Then
\begin{multline*}
(\forall k,l\in\N)(\forall x\in [a+h_1+h_2, b-h_1-h_2])\\
(x+kh_1-lh_2\in[a,b]\implies \abs{f(x+kh_1-lh_2)-f(x)}\le (k+l)\eps).
\end{multline*}
\end{lemma}
\begin{proof}
Let $x\in[a+h_1+h_2, b-h_1-h_2]$. Then the value of $f$ differs at most by
$\eps$ every time we move with a step $\pm h_1$ or $\pm h_2$, as long as we ensure that
all points lie in $[a,b]$. So it is sufficient to ensure that we reach
$x+kh_1-lh_2$ in at most $k+l$ such steps. We take the following steps: $x+h_1$,
$x+2h_1$, \dots, $x+k'h_1$, where $k'$ is the largest number such that
$x+k'h_1\le b$ (by the hypotheses, at least one step is taken). Then we move to
$x+k'h_1-h_2$, $x+k'h_1-2h_2$, \dots, $x + k'h_1 - l'h_2$, where $l'$ is the
largest number such that $x+k'h_1-l'h_2\ge a$ (again by the hypotheses, at least
one step is taken). Repeating this procedure, the coefficients of $h_1$ and $h_2$
increase until either the coefficient of $h_1$ equals $k$ or the coefficient of
$h_2$ equals $l$. In the first case, we have $x+ kh_1 -l'h_2$ $\in [a,b]$ with
$0\le l'\le l$. By hypothesis, also $x + kh_1 -l h_2$ $\in[a,b]$, so all the
remaining steps $x + kh_1 - (l'+1) h_2$, \dots, $x+ kh_1 + (l-1)h_2$ also lie in
$[a,b]$. The second case is similar.
\end{proof}

\begin{thm}
Let $f\in\Gen(\R)$. Let $\alpha\in\R^+\setminus\Q$ be an algebraic number.
If $f$ is periodic with periods $1$ and $\alpha$,
i.e., if $f(x+1)=f(x)$ and $f(x+\alpha)=f(x)$ hold as equalities in $\Gen(\R)$,
then $f$ is a generalized constant.
\end{thm}
\begin{proof}
Let $p\in\N$ and $K=[-R,R]\csub\R$.
By the corollary to thm.~\ref{Liouville}, we find $\forall\eps\in(0,1/\sqrt[p]2)$, some
$k_\eps$ and $l_\eps\in\N$, $l_\eps\le 1/\eps^p$ such that
\[
\eps^{Mp}\le\abs{k_\eps-\alpha l_\eps}\le 2\eps^p.
\]
In particular, $k_\eps\le \alpha l_\eps + 2\eps^p\le (\alpha + 1)/\eps^p$.\\
By the hypothesis, there exists $\eps_0>0$ such that for each $\eps\le\eps_0$,
$f_\eps$ is almost-periodic on $K$ with periods $1$, $\alpha$ and tolerance
$\eps^{(M+2)p}$. By lemma~\ref{almost-periodic}, we conclude that for $x\in\R$ with
$\abs{x}\le R-\alpha-1$,
\[\abs{f_\eps(x+k_\eps-\alpha l_\eps)-f_\eps(x)}\le (k_\eps + l_\eps)\eps^{(M+2)p}\,,\]
since $\abs{x+k_\eps-\alpha l_\eps}\le\abs{x} + 2\eps^p \le R$.
Let $h_\eps:=\abs{k_\eps-\alpha l_\eps}$. Now for each $x\in\R$ with
$\abs{x}\le R-\alpha - 2$,
we can find $\lambda_\eps\in\Z$ with $\abs{\lambda_\eps-x/h_\eps}\le 1$,
so also $\abs{\lambda_\eps h_\eps-x}\le 2\eps^p$. Then for each $j\in\Z$ with $\abs{j}\le\abs{\lambda_\eps}$,
$\abs{jh_\eps}\le R-\alpha - 1$, so
\begin{align*}
\abs{f_\eps(\lambda_\eps h_\eps) - f_\eps(0)}&\le
\sum_{j=1}^{\abs{\lambda_\eps}} \abs{f_\eps(\sigma_{\lambda_\eps}j
h_\eps)-f_\eps(\sigma_{\lambda_\eps}(j-1) h_\eps)}\\
&\le\abs{\lambda_\eps}(k_\eps+l_\eps)\eps^{(M+2)p}
\le \left(\frac{\abs{x}}{h_\eps}+1\right)(\alpha+2)\eps^{(M+1)p}
\le c\eps^p,
\end{align*}
for some constant $c$ only depending on $R$ (here $\sign_{\lambda_\eps}=\pm 1$ is the sign of
$\lambda_\eps$).\\
By the moderateness of $f'$, there exists $N\in\N$ (only depending on $R$)
such that
$\abs{f_\eps(x)- f_\eps(\lambda_\eps h_\eps)}\le\abs{x-\lambda_\eps h_\eps}\eps^{-N}\le
2\eps^{p-N}$ as soon as $\eps\le\eps_0$ (possibly with a smaller $\eps_0$).
So
\[
\sup_{\abs{x}\le R-\alpha-2}\abs{f_\eps(x)-f_\eps(0)}
\le c\eps^p + 2\eps^{p-N}, \quad \eps\le\eps_0\,.
\]
As $R$ and $p$ are arbitrary, it follows that $f$ is a generalized constant.
\end{proof}

\begin{cor}
Let $f\in\Gen(\R)$.
If $f$ is periodic with periods $h_1$ and $h_2$, with $h_1/h_2\in\R^+\setminus\Q$
algebraic, then $f$ is a generalized constant.
\end{cor}

{\bf Question}: suppose $f\in\Gen(\R)$ is periodic with periods $h_1$, $h_2\in\R^+$ and suppose that $h_1/h_2\notin\Q$. Is $f$ a generalized constant?\\
A generalization of the corollary to thm.~\ref{Liouville} for arbitrary $h\in\R^+\setminus\Q$ instead of
algebraic numbers is sufficient. Notice that an approximation as in Liouville's
approximation theorem, and hence also the corollary, holds for many
transcendent numbers as well~\cite{Besicovitch} (the exceptions form a set of
Hausdorff dimension 0).

We conclude with the nonstandard version.
\begin{lemma}\label{cnt}
Let $\Omega\subseteq\R^d$, $\Omega$ open.
For an internal map $f$: $\ster\Omega\to\ster\C$, the following are equivalent:
\begin{enumerate}
\item $(\forall x,y\in\ns(\ster\Omega))(x\approxeq y \implies f(x)\approxeq f(y))$
\item $(\forall K\csub\Omega)(\forall x\in\ster K)(\forall n\in\N)(\exists m\in\N)(\forall
y\in\ster K)(\abs{x-y}<\rho^m\implies \abs{f(x)-f(y)}<\rho^n)$
\end{enumerate}
In such case, we say that $f$ is $\approxeq$-continuous.
\end{lemma}
\begin{proof}
$\Rightarrow$: Let $K\csub\Omega$, $x\in\ster K$ and $n\in\N$.
By underspill, the internal set
\[
\{m\in\ster\N: (\forall y\in\ster K)(\abs{x-y}<\rho^m \implies
\abs{f(x)-f(y)}<\rho^n)\}
\]
contains some $m\in\N$.\\
$\Leftarrow$: Let $x$, $y$ $\in\ns(\ster\Omega)$.
Then there exists $K\csub\Omega$ such that $x,y\in \ster K$. If $x\approxeq
y$, we have $\abs{x-y}<\rho^m$, $\forall m\in\N$, so the hypothesis learns that
$(\forall n\in\N)(\abs{f(x)-f(y)}<\rho^n)$, i.e., $f(x)\approxeq f(y)$.
\end{proof}

\begin{thm}\label{ns-transinvar}
Let $f$: $\ster\R\to\ster\C$ be internal and $\approxeq$-continuous.
Let $\alpha\in\R^+\setminus\Q$ be an algebraic number.
Suppose that for $h\in\{1,\alpha\}$,
$f$ is almost-periodic up to iotas with period $h$, i.e.,
\[
(\forall x\in\ns(\ster\R))(f(x+h)\approxeq f(x)).
\]
Then $f$ is constant up to iotas, i.e.,
\[
(\forall x\in\ns(\ster\R))(f(x)\approxeq f(0)).
\]
\end{thm}
\begin{proof}
By transfer on the corollary to thm.~\ref{Liouville},
we find for each $n\in\N$ some $k$, $l$ $\in\ster\N$,
$l\le 1/\rho^n$, such that (still with the same standard $M$)
\[
\quad \rho^{nM} \le \abs{k-\alpha l}\le 2\rho^n\,.
\]
In particular, $k$, $l$ are moderate.\\
Let $R\in\R^+$. By overspill, the almost-periodicity up to iotas implies that
there exists some $\omega\in\ster\N\setminus\N$ such that for $h\in\{1, \alpha\}$
\[
(\forall x\in\ster\R)(\abs{x}\le R \implies \abs{f(x+h)-f(x)}\le \rho^\omega).
\]
By transfer on lemma~\ref{almost-periodic},
\begin{multline*}
(\forall k,l\in\ster\N)(\forall x\in \ster[-R+\alpha + 1, R-\alpha - 1])\\
(\abs{x+k-l\alpha}\le R\implies \abs{f(x+k-l\alpha)-f(x)}\le (k+l)\rho^\omega).
\end{multline*}
If $k$, $l$ are moderate and $k-l\alpha\approx 0$, we have in particular (as $R$ can be taken
arbitrarily large) that $k-l\alpha$ is an almost-period up to iotas for $f$ on the
whole of $\ns(\ster\R)$.\\
So for each $n\in\N$, $f$ has almost-periods up to iotas $h_n$ with
$\rho^{nM}<h_n<2\rho^n$.
Since $h_n\not\approxeq 0$, we find for each $x\in\ns(\ster\R)$ a moderate
$\lambda\in\ster\Z$ such that $\abs{\lambda h_n - x}<2\rho^n$.
Then
\[\abs{f(\lambda
h_n)-f(0)}\le\sum_{j=1}^{\abs{\lambda}}\abs{f(\sign_\lambda
jh_n)-f(\sign_\lambda(j-1)h_n)},\]
so $f(\lambda h_n)\approxeq f(0)$ (here $\sign_\lambda=\pm 1$ is the sign of
$\lambda$).
Since $f$ is internal and $\approxeq$-continuous, by the second characterization in
lemma~\ref{cnt}, this implies that $f(x)\approxeq f(0)$.
\end{proof}

\begin{cor}
Let $f\in {}^\rho{\mathcal E}(\R)$. If $f$ is periodic with periods $h_1$ and $h_2$, with
$h_1/h_2\in\R^+\setminus\Q$ algebraic, then $f$ is constant.
\end{cor}

\begin{cor}
Let $f\in {}^\rho{\mathcal E}(\R^d)$. If $f$ is invariant under standard translations,
i.e., $f(x+h)=f(x)$ holds in ${}^\rho{\mathcal E}(\R^d)$, $\forall h\in\R^d$, then $f$ is constant.
\end{cor}
\begin{proof}
Let $\tilde f\in\ster{\Cnt\infty}(\R^d)$ be a representative of $f$ and let $a=(a_1,\dots,a_d)\in\ns(\ster\R^d)$. Let $\tilde g(t)=\tilde f(a_1,\dots,a_{d-1},t)$, $t\in\ster\R$. Then $\tilde g\in\ster{\Cnt\infty}(\R)$ is a representative of an element $g\in {}^\rho{\mathcal E}(\R)$. For standard $h\in\R$ and $t\in\ns(\ster\R)$,
$\tilde g(t+h)=\tilde f(a_1,\dots,a_{d-1},t+h)\approxeq \tilde f(a_1,\dots,a_{d-1},t)=\tilde g(t)$, as $f$ is invariant under standard translations. This means that $g$ is invariant under standard translations, and is constant. In particular, $\tilde f(a)=\tilde g(a_d)\approxeq\tilde g(0)=\tilde f(a_1,\dots,a_{d-1}, 0)$. Similarly, we obtain
$\tilde f(a)\approxeq \tilde f(a_1,\dots,a_{d-1}, 0)\approxeq \tilde f(a_1,\dots,a_{d-2}, 0, 0)\approxeq\cdots\approxeq \tilde f(0)$. As $a\in\ns(\ster\R^d)$ arbitrary, $f$ is constant.
\end{proof}

\end{document}